\documentclass[12pt,a4paper,twoside]{article}                                % Artikel, 11pt-Schrift, DIN A4, etc.
\usepackage{graphicx}
\usepackage[latin1]{inputenc}                                               % Encoding "Manuskript" --> deutsche Umlaute und Sonderzeichen ("sz") werden erkannt
\usepackage[T1]{fontenc}                                                    % Encoding "Setzkasten" --> sog. Font-Encoding (?), standard
\usepackage{lmodern}                                                        % Latin Modern Font für korrekte Darstellung am Bildschirm (s. HU-Berlin Richtlinie)
\usepackage{upgreek,amsmath,amssymb,amstext}                         % ams-Pakete, d.h. weitere mathematische Symbole
\usepackage[amsmath,thmmarks]{ntheorem}                                     % Theoreme, Sätze, Beweise
\usepackage[OMLmathsfit]{isomath}                                           % kursive griechische Großbuchstaben, fett kursiv serifenlos etc.
\usepackage{accents}                                               % Komplexe Accents
\usepackage[total={15cm,23cm},top=3cm,centering]{geometry}            % eigenes Seitenlayout
\usepackage{parskip}                                                        % vertikaler Abstand zwischen Absätzen
\usepackage{url}
\usepackage[svgnames,x11names]{xcolor}                                       % Viele Farben im Dokument
\usepackage[vlined,ruled,linesnumbered]{algorithm2e}                                 % "
\usepackage{parskip}                                                    % paragraph setting: noindent + parskip (do NOT use \setlength{\parindent}{0em} + \setlength{\parskip}{0.6em})
\usepackage{siunitx}                                                   % Zahlen und Einheiten
\usepackage{setspace}

\newcommand{\bs}{\mathbold}                                               % fette kursive Symbole
\theoremstyle{plain}                        %
\newtheorem{theorAu}{Theorem}   %

\begin{document}

\title{\Large{Optimization-based smoothing algorithm for triangle meshes over arbitrarily shaped domains}}                              % title and author's names

\author{\large{D.~Aubram}\footnote{Correspondence to: Dr.-Ing.~Daniel Aubram, Chair of Soil Mechanics and Geotechnical Engineering, Technische Universit\"at Berlin (TU Berlin),
Secr.~TIB1-B7, Gustav-Meyer-Allee~25, D-13355~Berlin, Germany. Tel.:~+49~(0)30~31472349; fax:~+49~(0)30~31472343; E-mail:~daniel.aubram@tu-berlin.de}\\ \small{Chair of Soil Mechanics and Geotechnical Engineering, TU Berlin, Germany}}                         %
\date{}                                                                     %
\maketitle                                                                  %
\thispagestyle{plain}                                                       % first page has no header or footer

%********* Main Text ****************
\subsection*{Abstract}                        % narrowed abstract
\par{This paper describes a node relocation algorithm based on
nonlinear optimization which delivers excellent results for both
unstructured and structured plane triangle meshes over convex as
well as non-convex domains with high curvature. The local
optimization scheme is a damped Newton's method in which the
gradient and Hessian of the objective function are evaluated
exactly. The algorithm has been developed in order to continuously
rezone the mesh in arbitrary Lagrangian-Eulerian (ALE) methods for
large deformation penetration problems, but it is also suitable
for initial mesh improvement. Numerical examples highlight the
capabilities of the algorithm.

\textbf{Keywords:} mesh; triangle; smoothing; optimization; large
deformation; arbitrary La\-grang\-ian-Eulerian\par}

\section{Introduction}

In every mesh-based numerical method the convergence of the
solution algorithms and the accuracy of the solution results
depend on the quality of the mesh. Mesh improvement usually
becomes necessary, at least in postprocessing the originally
generated mesh. Mesh improvement is often initiated if a quality
measure drops below a certain value specified by the user.
Physical quality measures are employed in the adaptive numerical
methods for initial boundary value problems. Geometric quality
measures, including the size, aspect ratio, and skew of a mesh
element, can be evaluated independently of the physical solution
and usually at lower computational costs
\cite{Liu1994,Fie2000,Knu2001,Fry2000}.

The quality improvement of a mesh can be governed by quality
evolution and is done by repeated application of appropriate
tools. Smoothing is a tool intended to improve mesh quality by
node relocation. It represents a class of homeomorphic maps
between meshes which keep the connectivity of the original mesh
unchanged. Smoothing plays a crucial role in the arbitrary
Lagrangian-Eulerian (ALE) methods used for large deformation
problems with interfaces in computational solid and fluid dynamics
\cite{Hir1974,Hug1981,Ben1989,Naz2008,Sav2008a,Aub2010,Aub2013a};
see \cite{Ben1992,Don2004} for a review. ALE methods combine the
advantages of the purely Lagrangian and purely Eulerian
approaches. The computational mesh is not fixed but can move
independent of the material at an arbitrary velocity prescribed by
the smoothing scheme.

Since ALE methods must frequently relocate the mesh nodes when
advancing solution of the considered problem in time, an essential
requirement for the smoothing scheme is efficiency with respect to
computational costs. Another requirement closely connected with
efficiency is locality, that means to process only a set of
flagged nodes which may vary between the time steps. When using
local procedures, attention must be drawn to the strategy in order
to globally smooth the mesh. This holds for all local improvement
tools. Any improved mesh entity may deteriorate the quality of
neighboring entities. The third requirement imposed on a smoothing
algorithm is stability. A stable smoothing algorithm will not
distort a mesh any more than it is currently distorted
\cite{Ben1989}. For an algorithm to be reliable, this should be
independent of the domain's shape.

Automatic mesh smoothing procedures which are not governed by
quality evolution are called direct or heuristic smoothing
algorithms. Examples include Laplacian smoothing \cite{Fry2000},
smoothing by weighted averaging \cite{Aym2001}, and Giuliani's
method \cite{Giu1982}. These methods provide closed-form
expressions for the new node location which is supposed to smooth
the associated ball of elements sharing that node. Even though
these methods are computationally attractive, they cannot ensure
quality improvement for arbitrarily shaped domains. As will be
shown later, the heuristic smoothing algorithms proposed in
\cite{Aym2001} and \cite{Giu1982} fail on a non-convexly distorted
mesh.

Another class of algorithms is referred to as physically-based
smoothing. In these algorithms, physical properties are assigned
to the mesh entities and then a specific initial boundary value
problem is solved over a dummy time step in order to determine the
nodal displacements. Examples of physically-based smoothing
methods are reported in \cite{Loe1986,Naz2008}. The success of
such a procedure, however, is by pure chance. It cannot be ensured
that any mesh processed is not worse than before, i.e.~that the
smoothing scheme is stable.

The draw\-backs of heuris\-tic and physi\-cally-based smoo\-thing
techniques when dealing with non-convex meshes can be avoided if
the new node positions are determined through an optimization
process. In contrast to the other two approaches,
optimization-based smoothing algorithms are governed by geometric
quality evolution using an objective function whose minimum is
associated with a properly smoothed mesh, or a part of it. Early
references include \cite{Bar1970,Brb1982} which are concerned with
the global optimization of two-dimensional structured grids. One
of the first mesh smoothing algorithms that use principles of
local optimization is developed in \cite{Ken1986}. Several
refinements of the local approach and its generalization to
unstructured three-dimensional meshes are provided, for example,
in \cite{Dar1993,Zav1996}. Development continued up to the
present, with focuses on unstructured quadrilateral meshes
\cite{Jou1997}, unstructured triangle meshes \cite{Brs2000b},
structured quadrilateral meshes \cite{Knu2002}, and general
unstructured polyhedral meshes \cite{Dya2005a}. These algorithms
share the basic structure of all optimization procedures
\cite[sec.~1.5]{Sun2006} but differ in the methods to determine
the descent direction and step size, and particularly in the
objective function.

The remainder of this paper is concerned with the development of
an optimization-based smoo\-thing al\-go\-rithm for
two-di\-men\-sio\-nal tri\-angle me\-shes, which shall be referred
to as the OSMOT (Optimization-based SMOothing of Triangle meshes)
algorithm. The mesh can be an originally generated mesh but the
main objective of this research is to efficiently smooth distorted
meshes over non-convex domains arising in ALE finite element
simulations of penetration problems. Section~\ref{sec2} describes
the procedure to globally improve the mesh. The global algorithm
encloses local algorithms to smooth the boundary mesh and the
internal mesh which are outlined in Sections~\ref{sec3} and
\ref{sec4}, respectively. Extensions of the algorithm are
discussed in Section~\ref{sec7}. The numerical examples presented
in Section~\ref{sec5} highlight that the new algorithm delivers
excellent results for both unstructured and structured plane
triangle meshes over convex as well as non-convex domains with
high curvature. The paper closes with some concluding remarks in
Section~\ref{sec6}.

\section{Global algorithm}\label{sec2}

\subsection{General setup and initialization}\label{sec21}

Let $\mathcal{M}$ be a two-dimensional triangle mesh in the
Euclidian space $\mathcal{S}=\mathbb{R}^{2}$ and let
$\mathcal{N}(\mathcal{M})$ be the set of all nodes in the mesh.
The position vector of a node $P_0\in\mathcal{N}(\mathcal{M})$ is
given by $\bs{x}_0=(x_0,y_0)^{\mathrm{T}}\in\mathbb{R}^{2}$ with
respect to the canonical basis of $\mathbb{R}^{2}$. The
superscribed ${\mathrm{T}}$ denotes the transpose of a matrix.
Some frequently used geometric primitives of triangles are
compiled in Appendix~\ref{sec-app-tri-01}

The current procedure assumes that all nodes of the mesh are
allowed to be moved, except for the boundary nodes that
essentially define the shape of the meshed domain. The set of
internal nodes lying in the interior of the mesh is denoted by
$\mathcal{N}_{\mathrm{int}}\subset\mathcal{N}$. The non-movable
boundary nodes divide the boundary into a number of
$n_{\mathrm{bnd}}$ distinct sub-boundaries, and the set of all
movable nodes of the $j$-th sub-boundary is denoted by
$\mathcal{N}^j_{\partial}$, with
$j\in\{1,\ldots,n_{\mathrm{bnd}}\}$.

The algorithms intended to smooth the interior of a mesh generally
can not directly be applied to the boundary mesh. In most cases
the quality improvement of a distorted boundary mesh can be
achieved by simple heuristic procedures. Weighted averaging
\cite{Aym2001} is used here, whereas a new optimization-based
procedure is applied to smooth the internal mesh. These are local
algorithms in order to render the global improvement of the whole
mesh more effective.

The implemented local algorithms require additional topological
information. In particular, the local algorithm for internal nodes
works on the ball of elements associated with some node
$P_0\in\mathcal{N}$. A ball is the disjoint union
$\mathcal{B}(P_0)\overset{\mathrm{def}}{=}\bigcup_{n_{\mathrm{el}}}\!\triangle(P_0)$
of all $n_{\mathrm{el}}$ elements $\triangle$ in a mesh sharing
$P_0$, the vertex of the ball. Locally, the numbering of the nodes
in each triangle element of the ball is reordered such that the
signed area of the element (\ref{eq01}) is positive and the
location of the local node $0\in\triangle$ in
$\mathcal{S}=\mathbb{R}^{2}$ coincides with that of $P_0$. The
reordering of the local node numbers ensures that for each element
$\triangle\!\subset\mathcal{B}$ the vertex of the ball can be
addressed by $\bs{x}_0\in\mathbb{R}^{2}$, the position vector of
the local node $0\in\triangle$.

\subsection{Selection of the nodes to be moved}

Smoothing is initiated if at least one mesh element fails a
quality check. Stated loosely, a geometrically high quality mesh
is made up of more or less equal-sized elements with low
distortion. The two main groups of geometric quality measures are
accordingly referred to as size measures and shape measures. The
group of shape measures includes measures for the aspect ratio and
skew of an element \cite{Knu2001}.

For simplicial elements a size measure can be established by
taking the ratio of a reference radius $R_{\mathrm{ref}}$ and the
circumcircle radius $R$:
\begin{equation}\label{eq-ale-86}
Q_{1}\overset{\mathrm{def}}{=}\frac{R_{\mathrm{ref}}}{R}\;.
\end{equation}
However, $Q_1$ is a fair size measure only if the physical element
is almost regular, since $R$ can be finite even if element volume
is not (degenerate element).

A widely-used and versatile shape measure for simplicial elements
because it covers aspect ratio and skew is the normalized radius
ratio of the incircle and circumcircle
\cite{Liu1994,Fie2000,Brs2000b}:
\begin{equation}\label{eq-ale-79}
Q_2\overset{\mathrm{def}}{=}m\frac{r}{R}\qquad\in[0,1]\;.
\end{equation}
The normalization factor $m$ is the dimension of the simplex, with
$m=2$ (triangle) or $m=3$ (tetrahedron). A simplicial element is
equilateral if $Q_{2}=1$, and has zero volume if $Q_{2}=0$.

The geometric quality $Q_\triangle$ of each mesh element
$\triangle\!\in\mathcal{M}$ is compared with a minimal acceptable
quality $Q_{\mathrm{min}}$. The nodes of the elements that fail
the quality check are flagged. Hence, the set of flagged nodes
intended for relocation is given by
\begin{equation}
\mathcal{N}'\overset{\mathrm{def}}{=}\{P\in\mathcal{N}(\mathcal{M})\,|\,\mbox{$P\in\mathcal{N}(\triangle)$
and $\triangle\!\in\mathcal{M}$ and
$Q_\triangle<Q_{\mathrm{min}}$}\}\,,
\end{equation}
with $\mathcal{N}'\subset\mathcal{N}$. Computational costs of the
global algorithm can be considerably reduced by processing only
the flagged nodes by the local smoothing algorithms. Because the
total number of boundary nodes is not very large in
two-dimensional meshes, however, it would be adequate to relocate
all movable boundary nodes without any prior quality check.

\subsection{Global iteration}

The globally improved mesh is obtained by looping over the flagged
nodes of the mesh repeatedly. Hence, smoothing of the whole mesh
is achieved in an iterative fashion. Alg.~\ref{alg-ale-07}
provides the pseudocode of the entire procedure.

\IncMargin{2em}
\begin{algorithm}[!htbp]                      % enter the algorithm environment
\KwIn{triangle mesh $\mathcal{M}$, locations of the nodes $\mathcal{N}(\mathcal{M})$} %
\KwOut{smoothed mesh} %
\BlankLine %
initialize $i=0$, specify $i_{\mathrm{max}}$ and $Q_{\mathrm{min}}$\; %
specify set of movable nodes $\mathcal{N}^j_{\partial}$ for every
sub-boundary $j\in\{1,\ldots,n_{\mathrm{bnd}}\}$\; %
specify set of movable internal nodes $\mathcal{N}_{\mathrm{int}}$\; %
\ForEach{$P_0\in\mathcal{N}^j_{\partial}$ with
$j\in\{1,\ldots,n_{\mathrm{bnd}}\}$}{ %
determine neighboring nodes $P_1,P_2$\; %
}%
\ForEach{$P_0\in\mathcal{N}_{\mathrm{int}}$}{ %
determine ball of elements $\mathcal{B}(P_0)=\bigcup_{n_{\mathrm{el}}}\!\triangle(P_0)$\; %
}%
loop elements and evaluate element quality $Q_\triangle$\; %
\lIf{$Q_\triangle<Q_{\mathrm{min}}$}{flag nodes of element (set of all flagged nodes is $\mathcal{N}'$)}\;%
\While{global iteration step $i \leq i_{\mathrm{max}}$}{ %
\ForEach{$P_0\in\mathcal{N}^j_{\partial}$ with
$j\in\{1,\ldots,n_{\mathrm{bnd}}\}$}{ %
determine new location to smooth boundary mesh (Alg.~\ref{alg-ale-03})\;%
}%
\ForEach{$P_0\in(\mathcal{N}'\cap\mathcal{N}_{\mathrm{int}})$}{ %
smooth ball $\mathcal{B}(P_0)$ of internal mesh (Alg.~\ref{alg-ale-06})\; %
}%
$i\gets i+1$\;%
}%
\caption{Global mesh smoothing.}          % give the algorithm a caption
\label{alg-ale-07}                      % and a label for \ref{} commands later in the document
\end{algorithm}\DecMargin{2em}

\section{Local algorithm for boundary nodes}\label{sec3}

The heuristic algorithm of Aymone \emph{et al.}~\cite{Aym2001}
efficiently smoothes mesh boundaries. The boundary node
$P_0\in\mathcal{N}^j_{\partial}$ intended for relocation has two
neighbors, $P_1$ and $P_2$, which also belong to the boundary. To
relocate $P_0$, one assumes that the three points lie on a
sufficiently smooth curve
\begin{equation}\label{eq-ale-65}
\begin{aligned}
c: \,\, [-1,1] \ & \rightarrow \ \mathcal{S}=\mathbb{R}^2\\
\xi \ & \mapsto \ c(\xi)\,,\;\exists c^{-1}\,,
\end{aligned}
\end{equation}
with $c(-1)=P_1$, $c(0)=P_0$, and $c(1)=P_2$. The position vector
of a point $c(\xi)$ is
$\bs{c}(\xi)\overset{\mathrm{def}}{=}\bs{x}(c(\xi))$. Now the
boundary curve through $P_1,P_0,P_2$ considered in \cite{Aym2001}
is a polynomial of degree two such that $\bs{c}(\xi)$, with
$\xi\in[-1,1]$, has the exact representation
\begin{equation}\label{eq-ale-181}
\bs{c}(\xi)=\sum_{k=0}^{2}N_k(\xi)\,\bs{x}_k\;.
\end{equation}
Here $\bs{x}_k$, $k\in\{0,1,2\}$, are the position vectors of
$P_k$ in $\mathcal{S}$, and $N_k$ are quadratic interpolation
functions for $\bs{c}$ and $P_k$ having the particular form
\begin{equation}\label{eq-ale-69}
N_0(\xi)\overset{\mathrm{def}}{=}1-\xi^2\,,\qquad
N_1(\xi)\overset{\mathrm{def}}{=}\tfrac{1}{2}(\xi^2-\xi)\,,\qquad\mbox{and}\qquad
N_2(\xi)\overset{\mathrm{def}}{=}\tfrac{1}{2}(\xi^2+\xi)\,.
\end{equation}

A straightforward quality measure for the local boundary mesh
formed by $P_1,P_0,P_2$ is
\begin{equation}
Q_{\mathrm{bnd}}\overset{\mathrm{def}}{=}\frac{\min(\|\bs{x}_1-\bs{x}_0\|,\|\bs{x}_2-\bs{x}_0\|)}{\max(\|\bs{x}_1-\bs{x}_0\|,\|\bs{x}_2-\bs{x}_0\|)}\qquad\in[0,1]\,,
\end{equation}
with $\|\bs{x}_k-\bs{x}_l\|=\sqrt{(x_k-x_l)^2+(y_k-y_l)^2}$ and
$k,l\in\{0,1,2\}$. $Q_{\mathrm{bnd}}=1$ means that the location of
the node $P_0$ equalizes the distances (best quality). In
\cite{Aym2001}, weighted averaging is applied to determine a
natural coordinate $\xi'_0\in[-1,1]$ of $P_0$ that smoothes the
boundary curve. By using the distances $\|\bs{x}_1-\bs{x}_0\|$ and
$\|\bs{x}_2-\bs{x}_0\|$ as weights one arrives at
\begin{equation}
\xi'_0=\frac{\|\bs{x}_2-\bs{x}_0\|-\|\bs{x}_1-\bs{x}_0\|}{\|\bs{x}_1-\bs{x}_0\|+\|\bs{x}_2-\bs{x}_0\|}\,.
\end{equation}
The new position vector $\bs{x}'_0$ of $P_0$ that smoothes the
boundary curve can be obtained from (\ref{eq-ale-181}) by using
the coordinate $\xi=\xi'_0$, so that $\bs{x}'_0=\bs{c}(\xi'_0)$.
The procedure is summarized in Alg.~\ref{alg-ale-03}.

\IncMargin{2em}
\begin{algorithm}[!htbp]                      % enter the algorithm environment
\KwIn{neighboring nodes $P_1,P_2$ of every
$P_0\in\mathcal{N}^j_{\partial}$} %
\KwOut{smoothed position of $P_0\in\mathcal{N}^j_{\partial}$} %
\BlankLine %
read locations of nodes $\bs{x}_0=\bs{x}(P_0)$, $\bs{x}_1=\bs{x}(P_1)$, and $\bs{x}_2=\bs{x}(P_2)$\;%
compute distances $\|\bs{x}_1-\bs{x}_0\|$ and $\|\bs{x}_2-\bs{x}_0\|$\;%
natural coordinate to equalize distances is
$\xi'_0=\frac{\|\bs{x}_2-\bs{x}_0\|-\|\bs{x}_1-\bs{x}_0\|}{\|\bs{x}_1-\bs{x}_0\|+\|\bs{x}_2-\bs{x}_0\|}$\;%
location of $P_0$ smoothing the local boundary mesh is \newline $\bs{x}'_0=\sum_{k=0}^{2}N_k(\xi'_0)\,\bs{x}_k$, with $N_k$ given by (\ref{eq-ale-69})\;%
\caption{Local smoothing for boundary nodes.}          % give the algorithm a caption
\label{alg-ale-03}                      % and a label for \ref{} commands later in the document
\end{algorithm}\DecMargin{2em}

\section{Local optimization algorithm for internal nodes}\label{sec4}

\subsection{General remarks}\label{sec41}

Finding the best location of mesh nodes in terms of geometric
element quality constitutes an optimization problem which can be
solved using optimization theory \cite{Sun2006,Lue2008}. Let
$\mathcal{X}\subset\mathbb{R}^m$ be a feasible region,
$\bs{x}\in\mathcal{X}$, and $f:\mathcal{X}\rightarrow\mathbb{R}$ a
function. The general optimization problem can then be stated as
follows:
\begin{equation*}
\mbox{minimize}\;f(\bs{x})\qquad\mbox{subject to}\;\bs{x}\in
\mathcal{X}\,.
\end{equation*}
The function $f$ is called the objective function, and the
optimization problem is called unconstrained if
$\mathcal{X}=\mathbb{R}^m$. In the remainder of this paper, the
objective function is assumed to be twice continuously
differentiable in the Fr{\'e}chet-sense on $\mathcal{X}$, i.e.~of
class $C^2$ such that its gradient
$\boldsymbol{\nabla}\!f(\bs{x})\in\mathbb{R}^m$ and its Hessian
$\bs{H}_{\!f}(\bs{x})\in\mathbb{R}^{m\times m}$ at point
$\bs{x}\in\mathcal{X}$ do exist.

Determination of global minimizers is challenging. However, it is
usually sufficient to determine a local minimizer and to iterate
the global minimum. In this context the following first- and
second-order conditions are of fundamental importance. Proofs can
be found in \cite[sec.~1.4]{Sun2006}.

\begin{theorAu}\label{theo-ale-01}
(i) If $f:\mathcal{X}\rightarrow\mathbb{R}$ is continuously
differentiable on $\mathcal{X}\subset\mathbb{R}^m$ and
$\bs{x}'\in\arg\min_{\bs{x}\in\mathcal{X}}f(\bs{x})$ is a local
minimizer of $f$, then
\begin{equation*}
\boldsymbol{\nabla}\!f(\bs{x}')=\boldsymbol{0}\,.
\end{equation*}

(ii) If $\bs{x}'\in\mathcal{X}$ is a local minimizer of $f$, and
$f$ is of class $C^2$ on $\mathcal{X}$, then
$\boldsymbol{\nabla}\!f(\bs{x}')=\boldsymbol{0}$ and the Hessian
$\bs{H}_{\!f}(\bs{x}')$ is positive semidefinite,
i.e.~$\bs{x}^\mathrm{T}\bs{H}_{\!f}\,\bs{x}\geq0$ for every
$\bs{x}\in\mathbb{R}^m$ with $\bs{x}\neq\boldsymbol{0}$.

(iii) Let $f$ be $C^2$ on $\mathcal{X}$, $\bs{x}'\in\mathcal{X}$,
$\boldsymbol{\nabla}\!f(\bs{x}')=\boldsymbol{0}$, and let
$\bs{H}_{\!f}(\bs{x}')$ be positive definite such that
$\bs{x}^\mathrm{T}\bs{H}_{\!f}\,\bs{x}>0$ for every
$\bs{x}\in\mathbb{R}^m$ with $\bs{x}\neq\boldsymbol{0}$, then
$\bs{x}'$ is a strict local minimizer.
\end{theorAu}

Finding a local minimizer to solve the optimization problem
usually is an iterative procedure. Let $\mathcal{J}\in\mathbb{N}$
be an index set and $j, j+1\in\mathcal{J}$. For a given
$\bs{x}^{j}$, the iterative procedure takes the form
\begin{equation}\label{eq-ale-87}
\bs{x}^{j+1}=\bs{x}^j+\lambda^j\bs{d}^j\,,
\end{equation}
where $\bs{d}\in\mathbb{R}^m$ is a descent direction of $f$ at
$\bs{x}$ satisfying
$(\boldsymbol{\nabla}\!f(\bs{x}))^{\mathrm{T}}\bs{d}<0$. Once a
starting point $\bs{x}^{j=0}$, a step size $\lambda^{j=0}>0$, and
a tolerance $\varepsilon>0$ have been specified, a termination
criterion of the form
\begin{equation}
\|\boldsymbol{\nabla}\!f(\bs{x}^j)\|<\varepsilon
\end{equation}
is checked. If this criterion is met, then
$\bs{x}^j\approx\bs{x}'$ is an approximate minimizer of $f$. If
the criterion is not met, the descent direction $\bs{d}^j$
supposed to point to the minimum of the objective function is
determined by some method. Thereafter, a so-called line search is
carried out in order to determine the step size $\lambda^j$
satisfying
\begin{equation}
f(\bs{x}^j+\lambda^j\bs{d}^j)<f(\bs{x}^j)\,.
\end{equation}
An effective step size rule is highly desirable to ensure a
sufficient decrease in the objective function. The iterative
procedure continues with the repeated evaluation of the
termination criterion using
$\bs{x}^{j+1}=\bs{x}^j+\lambda^j\bs{d}^j$.

\subsection{Objective function}

The choice of an objective function is crucial to the success of
op\-ti\-mi\-za\-tion-based mesh smoothing. It must be composed of
geometric quality measures to ensure that the optimization is
governed by quality evolution. The class of local objective
functions for triangle meshes considered here takes the form
\cite{Brs2000b}
\begin{equation}\label{eq-ale-88}
W(\bs{x}_{0})\overset{\mathrm{def}}{=}\sum_{n_{\mathrm{el}}}w(\bs{x}_{0})\,,\qquad\mbox{with}\quad
w(\bs{x}_{0})\overset{\mathrm{def}}{=}\left(\frac{R(\bs{x}_{0})}{R_{\mathrm{ref}}}\right)^\beta\left(\frac{R(\bs{x}_{0})}{r(\bs{x}_{0})}\right)^\gamma\;.
\end{equation}
$n_{\mathrm{el}}$ is the number of triangles in the ball
$\mathcal{B}(P_0)=\bigcup_{n_{\mathrm{el}}}\triangle(P_0)$
associated with the internal node $P_0\in\mathcal{N}'$ whose
position vector is $\bs{x}_{0}$, $\beta$ and $\gamma$ are constant
positive weighting exponents, and $R_{\mathrm{ref}}>0$ is a
constant reference radius.

The class of local objective functions defined through
(\ref{eq-ale-88}) takes into account the size quality measure
(\ref{eq-ale-86}) and the shape quality measure (\ref{eq-ale-79}).
The weighting exponents $\beta$ and $\gamma$ control the
domination of the worst element. For example, if $\gamma$ is large
and $\beta$ is moderate, then the most distorted element
contributes more to the sum than a too large element or any of the
remaining elements. For the purpose of the present work, the
values $\beta=1.0$, $\gamma=3.0$, and $R_{\mathrm{ref}}=1.0$ have
been assigned to all elements in a mesh; see also \cite{Brs2000b}.

\subsection{Descent direction}

The first-order necessary condition (Theorem~\ref{theo-ale-01}(i))
in conjunction with (\ref{eq-ale-88}) defines a homogeneous system
of generally nonlinear algebraic equations,
$\boldsymbol{\nabla}W(\bs{x}_0')=\boldsymbol{0}$, whose solution
is $\bs{x}'_0\in\arg\min_{\bs{x}_0\in\mathcal{X}}W(\bs{x}_0)$, the
new location of the vertex of the ball $\mathcal{B}(P_0)$. The
solution can be approximated by Newton's method. Let $\bs{x}^j_0$
be a close-enough guess of the local minimizer. For
$W:\mathcal{X}\rightarrow\mathbb{R}$ being a $C^2$-function in the
neighborhood of $\bs{x}^j_0$, the linearization of
$\boldsymbol{\nabla}W(\bs{x}'_0)=\boldsymbol{0}$ about
$\bs{x}^j_0$ leads to
\begin{equation}
\boldsymbol{\nabla}W(\bs{x}_0^j)+\bs{H}_{\!W}(\bs{x}_0^j)\cdot
(\bs{x}_0'-\bs{x}_0^j)\approx\boldsymbol{0}\,.
\end{equation}
Provided that the gradient $\boldsymbol{\nabla}W$ is a
sufficiently smooth function, then any guess $\bs{x}_0^{j+1}$ for
which
$\boldsymbol{\nabla}W(\bs{x}_0^j)+\bs{H}_{\!W}(\bs{x}_0^j)\cdot(\bs{x}_0^{j+1}-\bs{x}_0^j)=\boldsymbol{0}$
is a better approximation than $\bs{x}_0^{j}$. If $\bs{H}_{\!W}$
is regular on $\mathcal{X}\subset\mathbb{R}^2$, this latter
condition results in the iterative scheme
\begin{equation}\label{eq-ale-73}
\bs{x}_0^{j+1}=\bs{x}_0^{j}-(\bs{H}_{\!W}^{-1}\boldsymbol{\nabla}W)(\bs{x}_0^j)\,,\qquad\mbox{and}\quad\lim_{j\rightarrow\infty}\bs{x}_0^{j+1}=\bs{x}_0'\,.
\end{equation}
By linearity,
\begin{equation}
\boldsymbol{\nabla}W(\bs{x}_0^j)=\sum_{n_{\mathrm{el}}}\boldsymbol{\nabla}w(\bs{x}_0^j)\qquad\mbox{and}\qquad\bs{H}_{\!W}(\bs{x}_0^j)=\sum_{n_{\mathrm{el}}}\bs{H}_{\!w}(\bs{x}_0^j)\,.
\end{equation}
Closed-form expressions for the components of the gradient
$\boldsymbol{\nabla}w(\bs{x}_0^j)$ and the Hessian
$\bs{H}_{\!w}(\bs{x}_0^j)$ are available through
(\ref{eq-ale-94})--(\ref{eq-ale-96}); see
Appendix~\ref{sec-app-smooth-01} for a straightforward
calculation.

By locality of Newton's method, the iterative scheme
(\ref{eq-ale-73}) converges only if the starting point
$\bs{x}_0^{j=0}$ is a close-enough guess of the solution. When the
starting point is far away from the solution it is not guaranteed
that the Hessian is invertible and positive definite at every
$\bs{x}_0^j\in\mathcal{X}$ and that Newton's direction,
$-(\bs{H}_{\!W}^{-1}\boldsymbol{\nabla}W)(\bs{x}_0^j)$, is indeed
a descent direction satisfying
\begin{equation}
(\boldsymbol{\nabla}W)^{\mathrm{T}}\bs{H}_{\!W}^{-1}\boldsymbol{\nabla}W>0\;.
\end{equation}
In these cases solution may diverge. Even if the starting point is
close to solution, the Hessian may still be non-positive definite
such that no strict local minimizer of the objective function
exists (cf.~Theorem~\ref{theo-ale-01}(iii)).

In order to ensure convergence at non-positive definite Hessians,
a modified Newton's method is employed. The particular approach
used in this work was suggested by Goldstein and
Price~\cite{Gol1967}. It substitutes the steepest descent
direction $-\boldsymbol{\nabla}W(\bs{x}_0^j)$ instead of
$-(\bs{H}_{\!W}^{-1}\boldsymbol{\nabla}W)(\bs{x}_0^j)$ for
$\bs{d}^j$ in (\ref{eq-ale-73}) whenever $\bs{H}_{\!W}$ is not
regular or positive definite at $\bs{x}_0^j$. The check for
positive definiteness is done by the angle criterion
\cite{Sun2006}. To this end, define
\begin{equation}\label{eq-ale-93}
\cos\theta^j\overset{\mathrm{def}}{=}-\frac{(\boldsymbol{\nabla}W(\bs{x}_0^j))^{\mathrm{T}}\bs{d}^j}{\|\boldsymbol{\nabla}W(\bs{x}_0^j)\|\|\bs{d}^j\|}
\end{equation}
at the $j$-th iteration. If $\cos\theta^j>0$ for
$\bs{d}^j=-(\bs{H}_{\!W}^{-1}\boldsymbol{\nabla}W)(\bs{x}_0^j)$,
then Newton's direction is indeed a descent direction
($\bs{H}_{\!W}$ is positive definite) and the iterative scheme
converges. If, on the other hand, $\cos\theta^j\leq0$ for
$\bs{d}^j=-(\bs{H}_{\!W}^{-1}\boldsymbol{\nabla}W)(\bs{x}_0^j)$,
then $\bs{d}^j=-\boldsymbol{\nabla}W(\bs{x}_0^j)$ is used as the
descent direction, satisfying $\cos\theta^j>0$ when substituted
into (\ref{eq-ale-93}) as long as $\bs{x}_0^j$ is not a minimizer
of $W$. This can be implemented as follows:
\begin{equation}\label{eq-ale-75}
\bs{d}^j\overset{\mathrm{def}}{=}\left\{%
\begin{array}{ll}
    -\boldsymbol{\nabla}W(\bs{x}_0^j), & \hbox{if $\det\bs{H}_{\!W}(\bs{x}_0^j)<\delta$ or if $\cos\theta^j<\eta$,} \\
    -(\bs{H}_{\!W}^{-1}\boldsymbol{\nabla}W)(\bs{x}_0^j), & \hbox{otherwise,} \\
\end{array}%
\right.
\end{equation}
where $\delta>0$ and $\eta>0$ are reasonable tolerances.

\subsection{Line search and step size rule}

It remains to determine the size of the steps with which the
optimization procedure approaches the local minimum of the
objective function. A too large step size may overshoot the
minimum, whereas a tiny step size would decelerate the overall
procedure. Therefore, it is a mandatory goal to let the programm
determine an appropriate size for every step by a line search.

Inexact line search is preferable from a computational viewpoint
provided that there is an effective step size rule which gives a
sufficient decrease in the objective function. One of such rules
is the widely-used Armijo rule \cite{Arm1966},
\begin{equation}\label{eq-ale-78}
W(\bs{x}_0^j+\lambda^j\bs{d}^j)-W(\bs{x}_0^j)\leq\tfrac{1}{2}\lambda^j(\boldsymbol{\nabla}W(\bs{x}_0^j))^{\mathrm{T}}\bs{d}^j\,,
\end{equation}
resulting in a so-called backtracking line search \cite{Sun2006}.
In a backtracking line search, for given $W(\bs{x}_0^j)$,
$\boldsymbol{\nabla}W(\bs{x}_0^j)$, $\bs{d}^j$, and
$\lambda^{j=0}=1.0$, the condition (\ref{eq-ale-78}) is checked.
If it is satisfied, then $\lambda^{j+1}=\lambda^j$ and
$\bs{x}_0^{j+1}=\bs{x}_0^j+\lambda^j\bs{d}^j$. If the condition is
not satisfied, the current guess of the minimizer is used in the
next iteration and the step size is bisected, that is,
\begin{equation}
\bs{x}_0^{j+1}=\bs{x}_0^j\qquad\mbox{and}\qquad\lambda^{j+1}=\frac{\lambda^j}{2}\,,
\end{equation}
respectively.

\subsection{Optimization procedure}

The entire local optimization procedure for internal nodes is
provided by Alg.~\ref{alg-ale-06}. The included tolerances have
been chosen to $\varepsilon=10^{-8}$, $\delta=10^{-6}$, and
$\eta=0.05$ for all the numerical examples presented in
Section~\ref{sec5}. Note that the algorithm is independent of the
specific objective function assigned to a ball of elements. It
might be attractive to implement alternative functions for which
evaluation of the gradient and Hessian is much cheaper or which
better reflect the user's needs.

%\begin{singlespace}
\IncMargin{2em}
\begin{algorithm}[!htbp]                      % enter the algorithm environment
\KwIn{ball $\mathcal{B}(P_0)$ associated with every
$P_0\in(\mathcal{N}'\cap\mathcal{N}_{\mathrm{int}})$} %
\KwOut{smoothed ball} %
\BlankLine %
specify tolerances $\varepsilon$, $\delta$, and $\eta$\; %
initialize $j=0$, $\bs{x}_0^{j=0}=\bs{x}(P_0)$, and $\lambda^{j=0}=1.0$\; %
\While{damped Newton iteration step $j \leq j_{\mathrm{max}}$}{ %
initialize $W(\bs{x}_0^j)\!=\!0$, $W(\bs{x}_0^j+\lambda^j\bs{d}^j)\!=\!0$, $\boldsymbol{\nabla}W(\bs{x}_0^j)\!=\!\boldsymbol{0}$, $\bs{H}_{\!W}(\bs{x}_0^j)\!=\!\boldsymbol{0}$\; %
\While{element in the ball $e \leq n_{\mathrm{el}}$}{ %
read and store locations of nodes $\bs{x}_0^j$, $\bs{x}_1$, and $\bs{x}_2$\;%
compute element objective function $w(\bs{x}_0^j)$\;%
compute $\boldsymbol{\nabla}w(\bs{x}_0^j)$ and $\bs{H}_{\!w}(\bs{x}_0^j)$ (\ref{sec-app-smooth-01})\;%
$W(\bs{x}_0^j)\gets W(\bs{x}_0^j)+w(\bs{x}_0^j)$\;%
$\boldsymbol{\nabla}W(\bs{x}_0^j)\gets\boldsymbol{\nabla}W(\bs{x}_0^j)+\boldsymbol{\nabla}w(\bs{x}_0^j)$\;%
$\bs{H}_{\!W}(\bs{x}_0^j)\gets\bs{H}_{\!W}(\bs{x}_0^j)+\bs{H}_{\!w}(\bs{x}_0^j)$\;%
}%
\eIf{$\|\boldsymbol{\nabla}W(\bs{x}_0^j)\|<\varepsilon$}{exit (location of $P_0$ is optimal)\;}{%
\textsc{Descent Direction:}\\%
\eIf{$\det\bs{H}_{\!W}(\bs{x}_0^j)<\delta$}{steepest descent $\bs{d}^j=-\boldsymbol{\nabla}W(\bs{x}_0^j)$\;}{%
Newton's direction $\bs{d}^j=-(\bs{H}_{\!W}^{-1}\boldsymbol{\nabla}W)(\bs{x}_0^j)$\;%
$\cos\theta^j=-(\boldsymbol{\nabla}W(\bs{x}_0^j))^{\mathrm{T}}\bs{d}^j/(\|\boldsymbol{\nabla}W(\bs{x}_0^j)\|\|\bs{d}^j\|)$\;%
\If{$\cos\theta^j<\eta$}{%
steepest descent $\bs{d}^j=-\boldsymbol{\nabla}W(\bs{x}_0^j)$\;%
}%
}%
\textsc{Line Search:}\\%
\While{element in the ball $e \leq n_{\mathrm{el}}$}{ %
compute element objective function $w(\bs{x}_0^j+\lambda^j\bs{d}^j)$\;%
$W(\bs{x}_0^j+\lambda^j\bs{d}^j)\gets W(\bs{x}_0^j+\lambda^j\bs{d}^j)+w(\bs{x}_0^j+\lambda^j\bs{d}^j)$\;%
}%
\eIf{$W(\bs{x}_0^j+\lambda^j\bs{d}^j)-W(\bs{x}_0^j)\leq\tfrac{1}{2}\lambda^j(\boldsymbol{\nabla}W(\bs{x}_0^j))^{\mathrm{T}}\bs{d}^j$}{update nodal location $\bs{x}_0^{j+1}=\bs{x}_0^j+\lambda^j\bs{d}^j$, whereas $\lambda^{j+1}=\lambda^j$\;}{%
update step size $\lambda^{j+1}=\tfrac{1}{2}\lambda^j$, whereas $\bs{x}_0^{j+1}=\bs{x}_0^j$\;%
}%
}%
}%
location of $P_0$ smoothing the ball is $\bs{x}'_{0}=\bs{x}^j_{0}$\;%
\caption{Local optimization-based smoothing for internal nodes.}          % give the algorithm a caption
\label{alg-ale-06}                      % and a label for \ref{} commands later in the document
\end{algorithm}\DecMargin{2em}
%\end{singlespace}

\section{Extensions to the current algorithm}\label{sec7}

The current algorithm can be naturally extended to adaptive
smoothing, to three dimensions, to higher-order elements, and to
surface meshes. The purpose of the following section is to discuss
these extensions.

The reference radius used in the objective function
(\ref{eq-ale-88}) is an attribute assigned to every element and
defines its maximum acceptable size in the mesh. Specifying
appropriate $R_{\mathrm{ref}}$ thus controls mesh grading during
the optimization process. If the value of the reference radius is
not specified by the user but \emph{a posteriori} by element
quality measures based on a numerical solution, then the
optimization-based algorithm would account for mesh grading,
leading to $r$-adaptive mesh improvement.

The proposed algorithm is based on a simplicial element type,
hence it should have a natural extension to three dimensions if
the triangles are replaced by tetrahedra. The local optimization
algorithm running over the internal nodes has a three dimensional
analog because all simplicial elements have a unique incircle and
a unique circumcircle, whose radii can be substituted into the
objective function (\ref{eq-ale-88}). In 3d, however, exact
evaluation of the gradient and Hessian of the related objective
function would yield awfully lengthy expressions. Moreover, for
the boundary nodes in 3d smoothing is much more complicated as one
will have to deal with surface triangulation connected to 3d
elements. A three-dimensional extension of the 2d algorithm
presented in Section~\ref{sec3} is proposed in \cite{Aym2001}, but
its suitability for simplicial meshes is not clear.

Extension to non-simplicial and/or higher-order element types with
midside nodes would generally require completely different
algorithms. However, a non-simplicial element can be divided into
simplices which can be processed by the current procedure. For
elements with midside nodes, a cheap but probably inadequate
approach would be to relocate only the corner nodes of an element
and to interpolate the midside nodes.

Smoothing of a surface triangle mesh by using the current
algorithm is non-trivial. One method to generate a surface mesh is
to regard the domain to be meshed as a parametric surface
\cite{Fry2000,Geo2004}, described by a two-dimensional parametric
domain $\mathcal{D}\subset\mathbb{R}^2$ and a smooth embedding
$\theta:\mathcal{D}\rightarrow\mathbb{R}^3$. Once the parametric
domain has been meshed, the map $\theta$ establishes the surface
mesh. The triangle mesh in the parametric space can be properly
smoothed by using the current algorithm, but the resulting quality
of the corresponding surface mesh is governed by $\theta$.

\section{Numerical examples}\label{sec5}

This section presents numerical examples highlighting the
applicability of the developed smoothing algorithm to different
types of meshes and mesh configurations. For reasons of
comparison, two additional smoothing algorithms for internal
meshes have been implemented. These heuristic algorithms may
replace Alg.~\ref{alg-ale-06} and likewise process the ball of
elements associated with a single internal node. The first is
based on weighted averaging and is used in the ALE method of
Aymone \emph{et al.}~\cite{Aym2001}. The second algorithm has been
developed by Giuliani \cite{Giu1982}. Its basic ingredient is an
objective function whose minimum yields closed-form expressions
for the coordinates of the internal node supposed to smooth the
ball.

\subsection{Patch tests}

\begin{figure}
\centering
\includegraphics{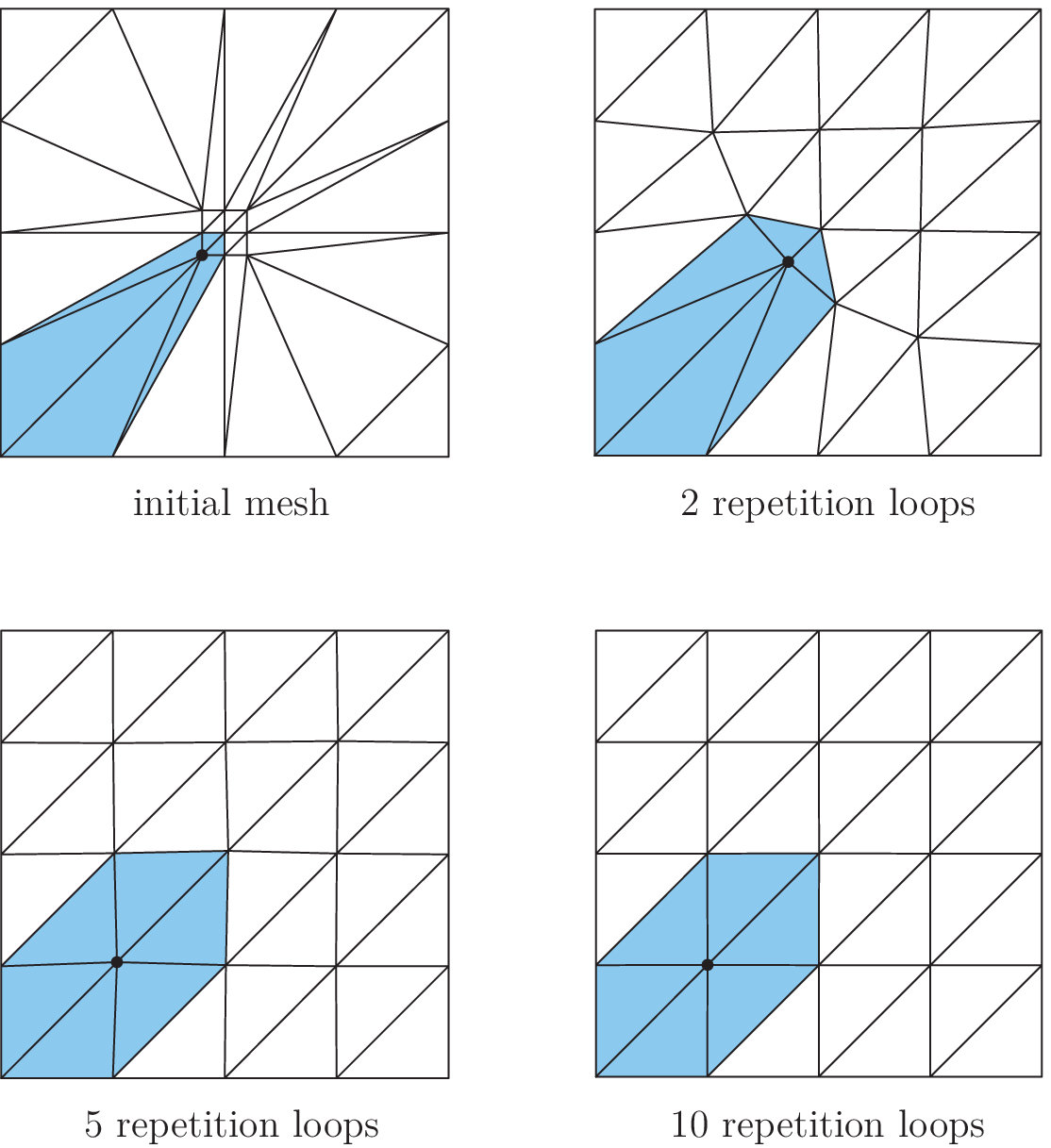}
\caption{Investigation of the number of repetition loops required
to globally improve the mesh. The blue zone indicates the ball of
elements sharing the lower left internal node.} \label{fig-num-01}
\end{figure}

The example shown in Fig.~\ref{fig-num-01} is a structured square
patch consisting of 32 triangle elements. The best quality of the
given mesh is obtained if the elements were arranged in a rising
diagonals triangle pattern. In the initial configuration, however,
elements are severely distorted. Merely the placement of the
boundary nodes is optimal. Due to locality of the implemented
smoothing algorithms, an acceptable mesh quality cannot be
achieved in one step but requires several repetition loops over
the balls of elements sharing a common internal node;
cf.~Alg.~\ref{alg-ale-07}. However, five to ten loops are
sufficient to produce an almost optimal mesh. This is independent
of the particular local smoothing algorithm used for the internal
mesh.

In the example, the quality improvement of the ball associated
with the lower left internal node (blue zone in
Fig.~\ref{fig-num-01}) lags behind the other after two iterations.
This is a consequence of the current strategy that globally
improves the mesh: all balls in the mesh are processed in a fixed
order in every repetition loop. It might be more effective to
process randomly picked groups of elements, but this has not been
implemented yet.

\begin{figure}
\centering
\includegraphics{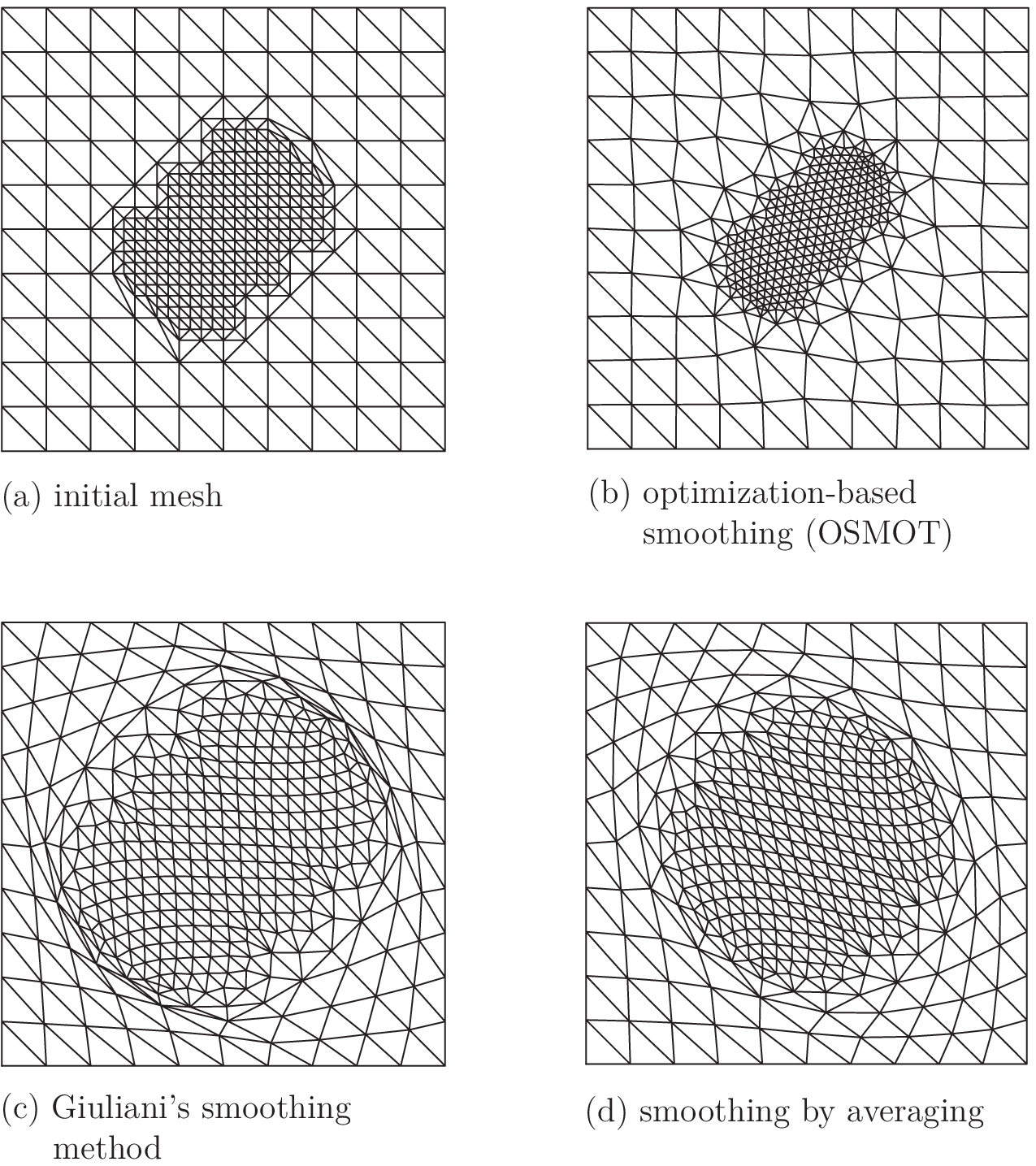}
\caption{Influence of the smoothing algorithm on mesh grading
after 1000 repetition loops.} \label{fig-num-02}
\end{figure}

The influence of the smoothing algorithm on mesh grading is
investigated in the second example. Graded or anisotropic meshes
made up of elements of prescribed size are often present in finite
element analysis, e.g.~when the computational model contains
regions of secondary interest. In these cases it is important to
preserve the prescribed element size. Fig.~\ref{fig-num-02}a shows
a structured triangle mesh zone with constant density interlaced
with a coarser structured mesh. The small interface zone is
unstructured and contains distorted elements, whereas the
structured parts of the mixed mesh are of best quality. An
appropriate smoothing algorithm hence would improve the interface
zone and would leave the structured zones unchanged.

It can be seen from Fig.~\ref{fig-num-02}b that after 1000 global
loops running over the internal nodes the optimization-based
algorithm (Alg.~\ref{alg-ale-06}) results in the mesh with the
highest quality. Structure is disturbed only slightly and the node
density distribution resp.~the size of elements is largely
preserved. In contrast to that, Giuliani's method \cite{Giu1982}
as well as smoothing by weighted averaging \cite{Aym2001} fail the
test (Fig.~\ref{fig-num-02}c and d, respectively). Both heuristic
procedures blow up the finer mesh zone, leading to a mesh with
equal-sized elements at repeated application. The quality of
elements at the interface deteriorated, Giuliani's method even
caused degenerate elements. The tendency to equalize the size of
elements is an undesirable feature which arises from the use of
averaged geometric measures in the governing equations of the
heuristic algorithms.

\subsection{Non-convexly distorted mesh}

A non-convexly distorted mesh is a mesh that contains stretched
and/or skewed elements in the vicinity of the indented boundary,
which probably has a high curvature. The automatic regularization
of such a mesh at fixed connectivity is very challenging. On the
other hand, problems associated with non-convexly distorted meshes
constitute important benchmark problems for the implemented
smoothing algorithms.

\begin{figure}
\centering
\includegraphics{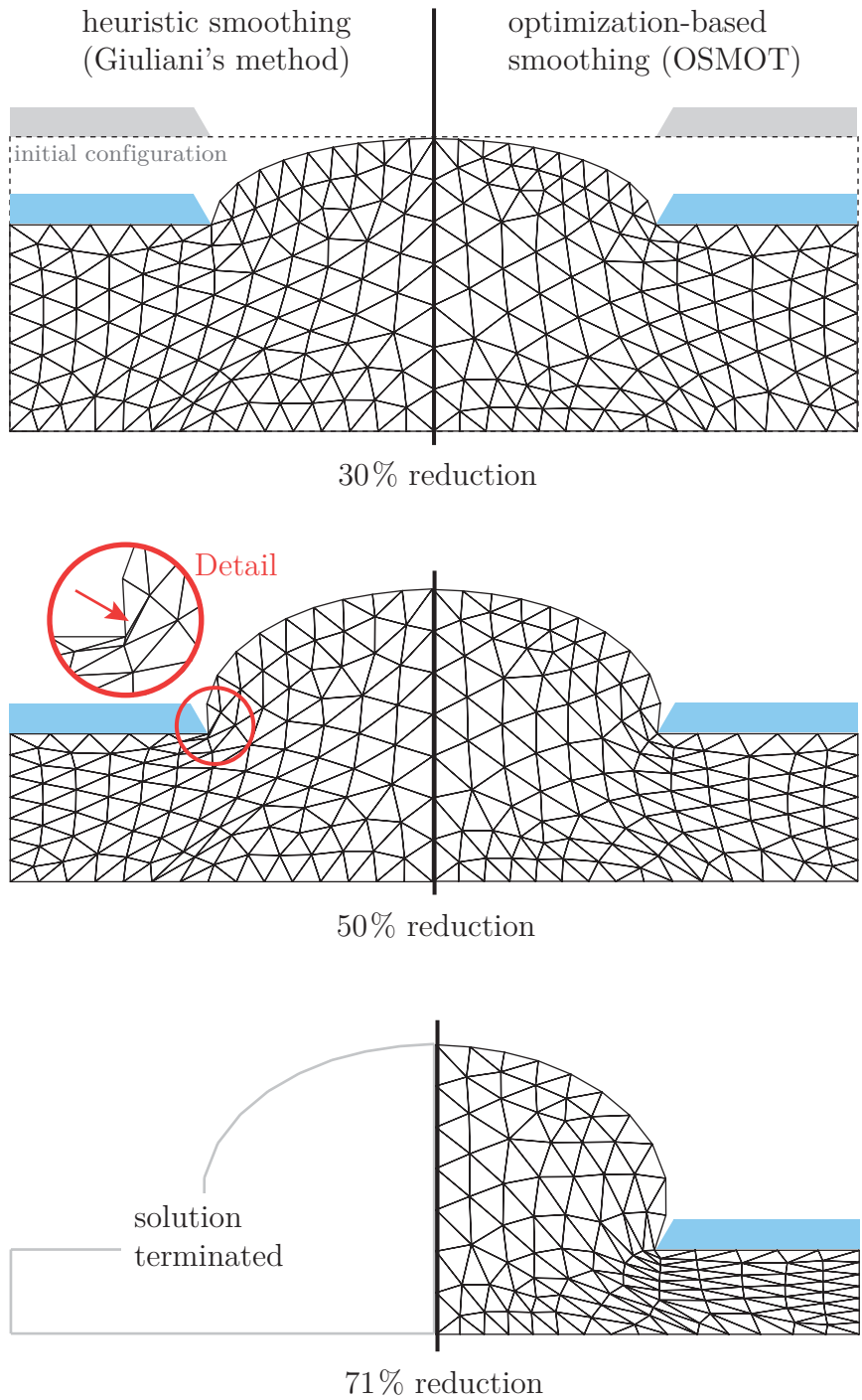}
\caption{Comparison of an heuristic smoothing method and the
developed op\-ti\-mi\-za\-tion-based algorithm when applied to the
numerical simulation of backward extrusion.} \label{fig-num-04}
\end{figure}

Backward extrusion is a common numerical example where non-convex
regions are created when large material deformation occurs. In
this initial boundary value problem a billet is loaded into a
heavy walled container and then a die is moved towards the billet,
so that the material is pushed through the die. Provided that the
die and the container are rigid and their surfaces are rough
respectively smooth, it suffices to discretize only the billet by
finite elements (Fig.~\ref{fig-num-04}). Nodes aligned with the
lower horizontal boundary are fixed in vertical direction, whereas
nodes at the wall of the container are fixed in horizontal
direction. The nodes located directly below the die are
horizontally fixed and will be displaced in vertical direction to
model the die moving downward.

Fig.~\ref{fig-num-04} above shows the edges of the undeformed
billet together with the deformed mesh at 30\,\% height reduction.
The left hand side shows the results of the calculation using a
heuristic scheme for mesh smoothing. Giuliani's method
\cite{Giu1982} for internal meshes has been employed, but
smoothing by weighted averaging would yield similar results. The
mesh on the right hand side results from the optimization-based
smoothing algorithm developed in this paper. In both calculations
the simple averaging procedure summarized in Alg.~\ref{alg-ale-03}
was chosen to smooth the boundary mesh.

The mesh quality of regions immediately under the die is
comparable at 30\,\% height reduction. Near the lower boundary,
the optimization-based algorithm produces a slightly smoother
mesh. At 50\,\% reduction, the heavy squeezing of elements around
the corner of the die cannot be avoided when using the heuristic
method. The area of one element even vanishes, which inhibits
convergence of the solution at continued extrusion. Compared to
the heuristic method, optimization-based smoothing achieves an
excellent mesh regularization. At 50\,\% height reduction, element
squeezing is moderate, even in the non-convexly distorted region
at the corner of the die. However, at continued extrusion the
fixed mesh connectivity associated with smoothing algorithms
limits gains of mesh quality. Calculation terminates at height
reductions of more than 71\,\%. Only a complete remeshing would
eliminate degenerate elements so as to continue solution.

\subsection{Penetration of a flat-ended pile}

The rigorous modeling of penetration is very challenging,
especially when the behavior of the penetrated material is highly
nonlinear. The final example is concerned with mesh smoothing
during the ALE simulation of the penetration of a flat-ended pile
into sand. It should be considered as an academic extreme example
highlighting the robustness of the new smoothing algorithm when
applied to large deformation problems involving indented material
boundaries. Details of the ALE method and the constitutive
equation used to model the mechanical behavior of sand can be
found in \cite{Aub2013a}.

\begin{figure}
\centering
\includegraphics{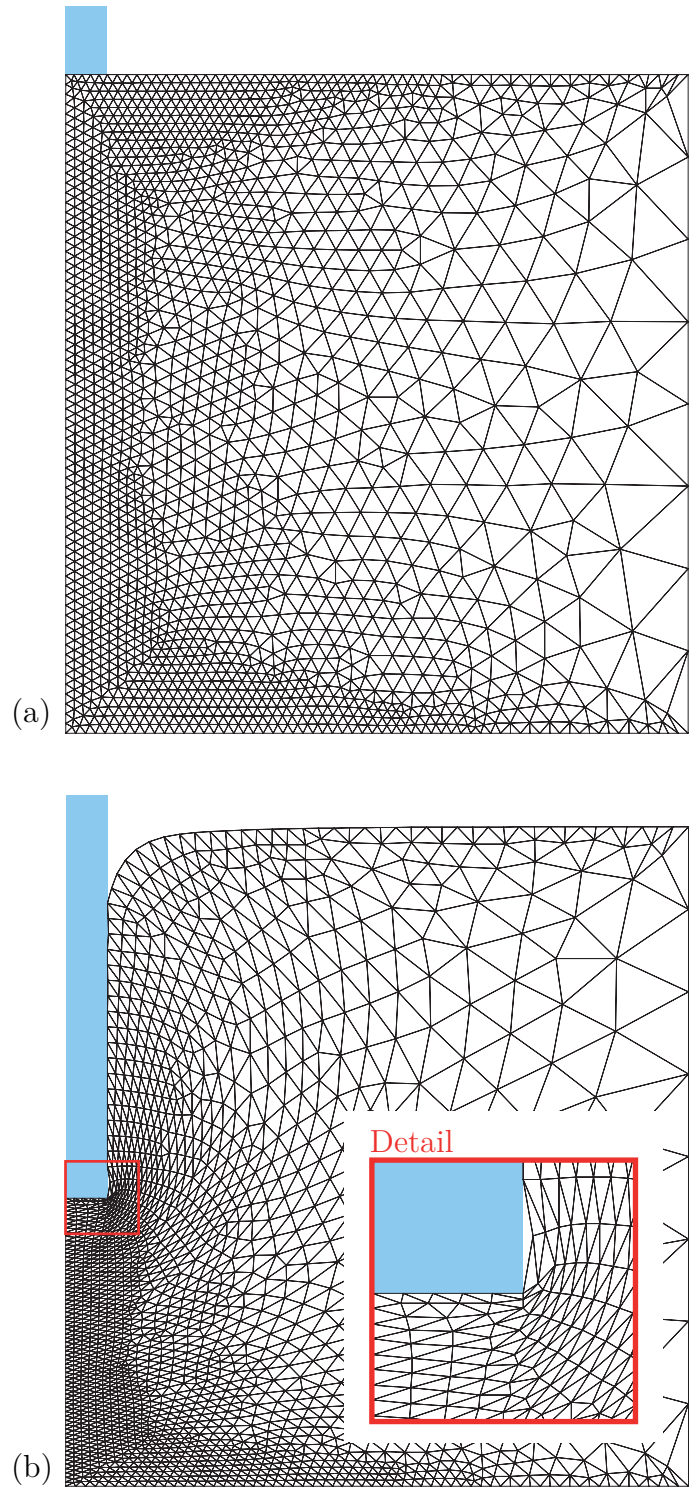}
\caption{ALE simulation of the penetration of a flat-ended pile
using the optimization-based algorithm. (a) Initial mesh, (b)
deformed and smoothed mesh at a relative penetration depth of
$z/D=4.5$.} \label{fig-num-20}
\end{figure}

The axisymmetric finite element model is depicted in
Fig.~\ref{fig-num-20}a. As penetration starts from the ground
surface, the initial configuration has a simple geometry. The pile
is assumed rigid, its shaft is assumed perfectly smooth, and the
pile base is perfectly rough (no sliding). The entire pile skin
and the ground surface are modeled as a contact pair using
straight segments for the sand surface and accounting for large
deformation of the interface. All nodes at the lower boundary of
the computational domain are fixed in vertical direction, and the
nodes at the vertical boundaries are fixed in radial direction.

For a relative penetration depth of $z/D=4.5$, where $D$ denotes
the diameter of the pile, the deformed and smoothed mesh is shown
in Fig.~\ref{fig-num-20}b. The initially rectangular computational
domain is severely deformed by indentation, resulting in a drastic
increase of the perimeter-to-area ratio. Elements at the elongated
boundary are stretched, i.e.~the density of nodes is reduced. On
the other hand, the local reduction in height of the domain below
the pile base comes along with squeezing of elements.
Optimization-based smoothing has indeed improved mesh quality in
this example. However, it could not completely avoid element
distortion because the mesh topology is fixed, that is, the mesh
underneath the pile cannot ``get out of the way''. Unless the
computational domain would be completely remeshed the numerical
model must contain a larger ``stockpile'' of less deformed mesh in
vertical direction in order to achieve a higher mesh quality.

\section{Conclusions}\label{sec6}

Heuristic smoothing algorithms, though they are simple and fast,
are inapplicable if the meshed domain becomes non-convex with high
curvature. Such situations may occur in initial mesh improvement
as well as in numerical simulations of large deformation problems.
In order to overcome these problems, an optimization-based
smoothing algorithm for triangle meshes has been developed in this
paper. The main objective was to continuously rezone the mesh in
an arbitrary Lagrangian-Eulerian method for penetration problems
\cite{Aub2013a}.

The smoothing algorithm operates iteratively on a local level and
distinguishes between boundary nodes and internal nodes. For
internal nodes, an optimization procedure has been developed which
processes the ball of elements enclosing a common node. It is
initiated if the quality measure based on the triangle's radius
ratio drops below a certain value specified by the user. A
globally smoothed mesh is pursued by loops repeatedly running over
the nodes of elements that fail the quality check. Numerical
examples show that the overall procedure, referred to as the OSMOT
(Optimization-based SMOothing of Triangle meshes) algorithm, is
efficient, extremely robust, and delivers excellent results for
both structured and unstructured triangle meshes over arbitrarily
shaped (i.e.~convex and non-convex) domains.

During penetration, an initially convex computational domain
necessarily becomes indented resp.~non-convex with high curvature
when material boundaries are explicitly resolved by element edges.
At drastic changes of the domain's shape due to large penetration
distances mesh quality improvement by smoothing can only be
achieved if there is a sufficiently large stockpile of less
deformed mesh. This would generally call for numerical models in
which the number of finite elements becoming additionally
necessary increases disproportionately with the desired
penetration distance. In this case, however, it is recommended to
completely remesh the computational domain or to try a different
modeling technique.

\section*{Acknowledgements}

The presented research work was carried out under the financial
support from the German Research Foundation (DFG), grant
SA~310/21-1, which is gratefully acknowledged.

\appendix


\setcounter{equation}{0}
\renewcommand{\theequation}{A\arabic{equation}}

\section{Geometric primitives of triangles}\label{sec-app-tri-01}

Consider a generic element
$\triangle\!\subset\mathcal{S}=\mathbb{R}^{2}$ of a plane triangle
mesh representing a 2-simplex. The local nodes $0,1,2\in\triangle$
occupy points $\bs{x}_k=(x_k,y_k)^{\mathrm{T}}\in\mathbb{R}^{2}$,
with $k\in\{0,1,2\}$. The local connectivity of the generic
triangle is predefined by choosing the first node $0$ and then
assigning the numbers of the other two nodes $1$ and $2$ in a
counter-clockwise manner such that the signed area $A$ given by
\begin{equation}\label{eq01}
2A=\det\left(%
\begin{array}{cc}
  x_1-x_0 & x_2-x_0 \\
  y_1-y_0 & y_2-y_0 \\
\end{array}%
\right)
\end{equation}
is positive. The lengths of the edges are readily available from
\begin{equation}
a\overset{\mathrm{def}}{=}\|\bs{x}_1-\bs{x}_0\|\,,\qquad
b\overset{\mathrm{def}}{=}\|\bs{x}_2-\bs{x}_1\|\,,\qquad\mbox{and}\qquad
c\overset{\mathrm{def}}{=}\|\bs{x}_0-\bs{x}_2\|\,,
\end{equation}
with $\|\bs{x}_k-\bs{x}_l\|=\sqrt{(x_k-x_l)^2+(y_k-y_l)^2}$ and
$k,l\in\{0,1,2\}$. The following relations for triangles are
well-known from undergraduate texts in geometry \cite{Bro2007}:
\begin{align}
\mbox{semiperimeter:}\qquad & s=\tfrac{1}{2}(a+b+c)\,,\label{eq-ale-94}\\
\mbox{Heron's formula:}\qquad & |A|=\sqrt{s(s-a)(s-b)(s-c)}\,,\\
\mbox{incircle radius:}\qquad & r=\dfrac{|A|}{s}\,,\label{eq-ale-95}\\
\mbox{circumcircle radius:}\qquad &
R=\dfrac{abc}{4|A|}\,.\label{eq-ale-96}
\end{align}
The functional dependence of $s,A,r,$ and $R$ on
$(\bs{x}_0,\bs{x}_1,\bs{x}_2)$ is being understood.

\setcounter{equation}{0}
\renewcommand{\theequation}{B\arabic{equation}}

\section{Gradient and Hessian for mesh
optimization}\label{sec-app-smooth-01}

\subsection{General remarks}

The optimization-based iterative local mesh smoothing algorithm
for internal nodes (Alg.~\ref{alg-ale-06}) requires frequent
evaluation of the gradient and Hessian of the element objective
function. By assuming that the numbering of the local nodes is in
accordance with Section~\ref{sec21} and assuming that the
locations $\bs{x}_{1},\bs{x}_{2}$ of the local triangle nodes
$1,2\in\triangle$ are constant during the iteration process, the
currently implemented objective function for a triangle in the
$j$-th iteration takes the form
\begin{equation}
w(\bs{x}^j_0)=\frac{R}{R_{\mathrm{ref}}}\left(\frac{R}{r}\right)^3=\frac{(abc)^4s^3}{4^4R_{\mathrm{ref}}A^7}\;,
\end{equation}
where (\ref{eq-ale-95}) and (\ref{eq-ale-96}) have been used, and
the functional dependence of $a,b,c,s,$ $A,r,R$ on
$(\bs{x}^j_0,\bs{x}_1,\bs{x}_2)$ is being understood. By dropping
the superscribed $j$ indicating the iteration step in what
follows, the gradient and Hessian of $w(\bs{x}_0)$ in
$\mathbb{R}^2$ are the component matrices given by
\begin{equation}
\renewcommand{\arraystretch}{1.5}
\boldsymbol{\nabla}w(\bs{x}_0)=\left(%
\begin{array}{c}
  \frac{\partial w}{\partial x_0} \\
  \frac{\partial w}{\partial y_0} \\
\end{array}%
\right)\qquad\mbox{and}\qquad\bs{H}_{\!w}(\bs{x}_0)=\left(%
\begin{array}{ccc}
  \frac{\partial^2 w}{\partial x^2_0} & \frac{\partial^2 w}{\partial x_0 \partial y_0}\\
  \frac{\partial^2 w}{\partial y_0 \partial x_0} & \frac{\partial^2 w}{\partial y^2_0} \\
\end{array}%
\right)\,,
\end{equation}
respectively. The geometric primitives of triangles provided in
\ref{sec-app-tri-01} enable the straightforward calculation of the
components of $\boldsymbol{\nabla}w(\bs{x}_0)$ and
$\bs{H}_{\!w}(\bs{x}_0)$, see below. Note that the derivatives of
$b=\sqrt{(x_2-x_1)^2+(y_2-y_1)^2}$ with respect to $\bs{x}_0$,
that is, $\frac{\partial b}{\partial x_0},\frac{\partial
b}{\partial y_0},\frac{\partial^2 b}{\partial x_0\partial y_0}$,
etc., identically vanish.

\subsection{First derivatives of objective function}

\begin{align}
\frac{\partial w}{\partial x_0}=\; &
\frac{1}{4^4R_{\mathrm{ref}}}\left(\frac{s^3}{A^7}\frac{\partial}{\partial
x_0}(abc)^4+\frac{(abc)^4}{A^7}\frac{\partial}{\partial
x_0}s^3+(abc)^4s^3\frac{\partial}{\partial
x_0}A^{-7}\right)\\
\frac{\partial w}{\partial y_0}=\; &
\frac{1}{4^4R_{\mathrm{ref}}}\left(\frac{s^3}{A^7}\frac{\partial}{\partial
y_0}(abc)^4+\frac{(abc)^4}{A^7}\frac{\partial}{\partial
y_0}s^3+(abc)^4s^3\frac{\partial}{\partial y_0}A^{-7}\right)
\end{align}

\subsubsection{Extensions}

{\allowdisplaybreaks
\begin{align}
\frac{\partial}{\partial x_0}(abc)^4=\; & 4b(abc)^3(cC_{1x}+aC_{2x})\\
\frac{\partial}{\partial y_0}(abc)^4=\; & 4b(abc)^3(cC_{1y}+aC_{2y})\\
\frac{\partial}{\partial x_0}s^3=\; & 3s^2C_{3x}\\
\frac{\partial}{\partial y_0}s^3=\; & 3s^2C_{3y}\\
\frac{\partial}{\partial
x_0}A^{-7}=\; & -\frac{7}{2A^9}\big[(2s-b)(s-a)(s-c)C_{3x}\notag\\
{} & +s(s-b)(s-c)C_{4x}+s(s-a)(s-b)C_{5x}\big]\\
\frac{\partial}{\partial
y_0}A^{-7}=\; & -\frac{7}{2A^9}\big[(2s-b)(s-a)(s-c)C_{3y}\notag\\
{} & +s(s-b)(s-c)C_{4y}+s(s-a)(s-b)C_{5y}\big]\\
\frac{\partial A}{\partial x_0}=\; &
-\frac{A^8}{7}\frac{\partial}{\partial x_0}A^{-7}\\
\frac{\partial A}{\partial y_0}=\; &
-\frac{A^8}{7}\frac{\partial}{\partial y_0}A^{-7}
\end{align}
}

\subsubsection{Abbreviations}

{\allowdisplaybreaks
\begin{align}
C_{1x}\overset{\mathrm{def}}{=}\; & \frac{\partial a}{\partial
x_0}=\frac{\partial}{\partial
x_0}\left((x_1-x_0)^2+(y_1-y_0)^2\right)^{\frac{1}{2}}=-\frac{x_1-x_0}{a}\\
C_{1y}\overset{\mathrm{def}}{=}\; & \frac{\partial a}{\partial
y_0}=\frac{\partial}{\partial
y_0}\left((x_1-x_0)^2+(y_1-y_0)^2\right)^{\frac{1}{2}}=-\frac{y_1-y_0}{a}\\
C_{2x}\overset{\mathrm{def}}{=}\; & \frac{\partial c}{\partial
x_0}=\frac{\partial}{\partial
x_0}\left((x_0-x_2)^2+(y_0-y_2)^2\right)^{\frac{1}{2}}=\frac{x_0-x_2}{c}\\
C_{2y}\overset{\mathrm{def}}{=}\; & \frac{\partial c}{\partial
y_0}=\frac{\partial}{\partial
y_0}\left((x_0-x_2)^2+(y_0-y_2)^2\right)^{\frac{1}{2}}=\frac{y_0-y_2}{c}\\
C_{3x}\overset{\mathrm{def}}{=}\; & \frac{\partial s}{\partial
x_0}=\frac{1}{2}\left(\frac{\partial a}{\partial
x_0}+\frac{\partial b}{\partial x_0}+\frac{\partial c}{\partial
x_0}\right)=\frac{1}{2}(C_{1x}+C_{2x})\\
C_{3y}\overset{\mathrm{def}}{=}\; & \frac{\partial s}{\partial
y_0}=\frac{1}{2}\left(\frac{\partial a}{\partial
y_0}+\frac{\partial b}{\partial y_0}+\frac{\partial c}{\partial
y_0}\right)=\frac{1}{2}(C_{1y}+C_{2y})\\
C_{4x}\overset{\mathrm{def}}{=}\; & \frac{\partial (s-a)}{\partial
x_0}=\frac{\partial s}{\partial x_0}-\frac{\partial a}{\partial
x_0}=C_{3x}-C_{1x}\\
C_{4y}\overset{\mathrm{def}}{=}\; & \frac{\partial (s-a)}{\partial
y_0}=\frac{\partial s}{\partial y_0}-\frac{\partial a}{\partial
y_0}=C_{3y}-C_{1y}\\
C_{5x}\overset{\mathrm{def}}{=}\; & \frac{\partial (s-c)}{\partial
x_0}=\frac{\partial s}{\partial x_0}-\frac{\partial c}{\partial
x_0}=C_{3x}-C_{2x}\\
C_{5y}\overset{\mathrm{def}}{=}\; & \frac{\partial (s-c)}{\partial
y_0}=\frac{\partial s}{\partial y_0}-\frac{\partial c}{\partial
y_0}=C_{3y}-C_{2y}
\end{align}
}

\subsection{Second derivatives of objective function}

{\allowdisplaybreaks
\begin{align}
\frac{\partial^2 w}{\partial x_0^2}=\; &
\frac{1}{4^4R_{\mathrm{ref}}}\left(\frac{s^3}{A^7}\frac{\partial^2}{\partial
x_0^2}(abc)^4+\frac{(abc)^4}{A^7}\frac{\partial^2}{\partial
x_0^2}s^3+s^3(abc)^4\frac{\partial^2}{\partial
x_0^2}A^{-7}\right.\notag\\
{} & +\frac{2}{A^7}\frac{\partial}{\partial
x_0}(abc)^4\frac{\partial}{\partial
x_0}s^3+2s^3\frac{\partial}{\partial
x_0}(abc)^4\frac{\partial}{\partial
x_0}A^{-7}\notag\\
{} & \left.+2(abc)^4\frac{\partial}{\partial
x_0}s^3\frac{\partial}{\partial x_0}A^{-7}\right)\\
\frac{\partial^2 w}{\partial y_0^2}=\; &
\frac{1}{4^4R_{\mathrm{ref}}}\left(\frac{s^3}{A^7}\frac{\partial^2}{\partial
y_0^2}(abc)^4+\frac{(abc)^4}{A^7}\frac{\partial^2}{\partial
y_0^2}s^3+s^3(abc)^4\frac{\partial^2}{\partial
y_0^2}A^{-7}\right.\notag\\
{} & +\frac{2}{A^7}\frac{\partial}{\partial
y_0}(abc)^4\frac{\partial}{\partial
y_0}s^3+2s^3\frac{\partial}{\partial
y_0}(abc)^4\frac{\partial}{\partial
y_0}A^{-7}\notag\\
{} & \left.+2(abc)^4\frac{\partial}{\partial
y_0}s^3\frac{\partial}{\partial y_0}A^{-7}\right)\\
\frac{\partial^2 w}{\partial x_0\partial y_0}=\; &
\frac{1}{4^4}\left(\frac{s^3}{A^7}\frac{\partial^2}{\partial
x_0\partial y_0}(abc)^4+\frac{1}{A^7}\frac{\partial}{\partial
x_0}(abc)^4\frac{\partial}{\partial
y_0}s^3\right.\notag\\
{} & +s^3\frac{\partial}{\partial
x_0}(abc)^4\frac{\partial}{\partial
y_0}A^{-7}+\frac{1}{A^7}\frac{\partial}{\partial
y_0}(abc)^4\frac{\partial}{\partial
x_0}s^3\notag\\
{} & +\frac{(abc)^4}{A^7}\frac{\partial^2}{\partial x_0\partial
y_0}s^3+(abc)^4\frac{\partial}{\partial
x_0}s^3\frac{\partial}{\partial
y_0}A^{-7}\notag\\
{} & +s^3\frac{\partial}{\partial
y_0}(abc)^4\frac{\partial}{\partial
x_0}A^{-7}+(abc)^4\frac{\partial}{\partial
y_0}s^3\frac{\partial}{\partial
x_0}A^{-7}\notag\\
{} & \left.+s^3(abc)^4\frac{\partial^2}{\partial x_0\partial
y_0}A^{-7}\right)=\frac{\partial^2 w}{\partial y_0\partial x_0}
\end{align}
}

\subsubsection{Extensions}

{\allowdisplaybreaks
\begin{align}
\frac{\partial^2}{\partial^2x_0}(abc)^4=\; & 12b(abc)^2(cC_{1x}+aC_{2x})^2\notag\\
{} & +4b(abc)^3(cD_{1x}+aD_{2x}+2C_{1x}C_{2x})\\
\frac{\partial^2}{\partial^2y_0}(abc)^4=\; & 12b(abc)^2(cC_{1y}+aC_{2y})^2\notag\\
{} &
+4b(abc)^3(cD_{1y}+aD_{2y}+2C_{1y}C_{2y})\\
\frac{\partial^2}{\partial x_0\partial
y_0}(abc)^4=\; & 12b(abc)^2(cC_{1x}+aC_{2x})(cC_{1y}+aC_{2y})\notag\\
{} &
+4b(abc)^3(cE_1+C_{1x}C_{2y}+C_{1y}C_{2x}+aE_2)\\
\frac{\partial^2}{\partial x_0^2}s^3=\; & 6sC_{3x}^2+3s^2D_{3x}\\
\frac{\partial^2}{\partial y_0^2}s^3=\; & 6sC_{3y}^2+3s^2D_{3y}\\
\frac{\partial^2}{\partial x_0\partial y_0}s^3=\; &
6sC_{3x}C_{3y}+3s^2E_3\\
\frac{\partial^2}{\partial x_0^2}A^{-7}=\; &
\frac{63}{A^{9}}\left(\frac{\partial A}{\partial
x_0}\right)^2\notag\\
{} & -\frac{7}{2A^9}\big[(s-a)\{(s-b)(sD_{5x}+(s-c)D_{3x}+2C_{3x}C_{5x})\notag\\
{} & +(s-c)(sD_{3x}+2C_{3x}^2)+2sC_{3x}C_{5x}\}\notag\\
{} & +(s-b)\{(s-c)(sD_{4x}+2C_{3x}C_{4x})+2sC_{4x}C_{5x}\}\notag\\
{} & +2s(s-c)C_{3x}C_{4x}\big]\\
\frac{\partial^2 }{\partial y_0^2}A^{-7}=\; &
\frac{63}{A^{9}}\left(\frac{\partial A}{\partial
y_0}\right)^2\notag\\
{} & -\frac{7}{2A^9}\big[(s-a)\{(s-b)(sD_{5y}+(s-c)D_{3y}+2C_{3y}C_{5y})\notag\\
{} & +(s-c)(sD_{3y}+2C_{3y}^2)+2sC_{3y}C_{5y}\}\notag\\
{} &
+(s-b)\{(s-c)(sD_{4y}+2C_{3y}C_{4y})+2sC_{4y}C_{5y}\}\notag\\
{} & +2s(s-c)C_{3y}C_{4y}\big]\\
\frac{\partial^2}{\partial x_0\partial y_0}A^{-7}=\; &
\frac{63}{A^{9}}\frac{\partial A}{\partial y_0}\frac{\partial
A}{\partial
x_0}\notag\\
{} & -\frac{7}{2A^9}\big[(s-a)\{(s-b)(C_{3x}C_{5y}\!+\!(s-c)E_3+C_{5x}C_{3y}+sE_5)\notag\\
{} & +(s-c)(2C_{3x}C_{3y}+sE_3)+s(C_{5x}C_{3y}+C_{3x}C_{5y})\}\notag\\
{} & +s(s-c)(C_{3x}C_{4y}+C_{4x}C_{3y})+(s-b)\{s(C_{4x}C_{5y}+C_{5x}C_{4y})\notag\\
{} & +(s-c)(C_{3x}C_{4y}+C_{4x}C_{3y}+sE_4)\}\big]
\end{align}
}

\subsubsection{Abbreviations}

{\allowdisplaybreaks
\begin{align}
D_{1x}\overset{\mathrm{def}}{=}\; & \frac{\partial^2 a}{\partial
x_0^2}=\frac{\partial}{\partial
x_0}\left(-\frac{x_1-x_0}{a}\right)=\frac{x_1-x_0}{a^2}C_{1x}+\frac{1}{a}\\
D_{1y}\overset{\mathrm{def}}{=}\; & \frac{\partial^2 a}{\partial
y_0^2}=\frac{\partial}{\partial
y_0}\left(-\frac{y_1-y_0}{a}\right)=\frac{y_1-y_0}{a^2}C_{1y}+\frac{1}{a}\\
D_{2x}\overset{\mathrm{def}}{=}\; & \frac{\partial^2 c}{\partial
x_0^2}=-\frac{x_0-x_2}{c^2}C_{2x}+\frac{1}{c}\\
D_{2y}\overset{\mathrm{def}}{=}\; & \frac{\partial^2 c}{\partial
y_0^2}=-\frac{y_0-y_2}{c^2}C_{2y}+\frac{1}{c}\\
D_{3x}\overset{\mathrm{def}}{=}\; & \frac{\partial^2 s}{\partial
x_0^2}=\frac{\partial}{\partial x_0}\left(\frac{\partial}{\partial
x_0}\frac{a+b+c}{2}\right)=\frac{1}{2}(D_{1x}+D_{2x})\\
D_{3y}\overset{\mathrm{def}}{=}\; & \frac{\partial^2 s}{\partial
y_0^2}=\frac{\partial}{\partial y_0}\left(\frac{\partial}{\partial
y_0}\frac{a+b+c}{2}\right)=\frac{1}{2}(D_{1y}+D_{2y})\\
D_{4x}\overset{\mathrm{def}}{=}\; & \frac{\partial^2
(s-a)}{\partial x_0^2}=\frac{\partial^2 s}{\partial
x_0^2}-\frac{\partial^2 a}{\partial x_0^2}=D_{3x}-D_{1x}\\
D_{4y}\overset{\mathrm{def}}{=}\; & \frac{\partial^2
(s-a)}{\partial y_0^2}=\frac{\partial^2 s}{\partial
y_0^2}-\frac{\partial^2 a}{\partial y_0^2}=D_{3y}-D_{1y}\\
D_{5x}\overset{\mathrm{def}}{=}\; & \frac{\partial^2
(s-c)}{\partial x_0^2}=\frac{\partial^2 s}{\partial
x_0^2}-\frac{\partial^2 c}{\partial x_0^2}=D_{3x}-D_{2x}\\
D_{5y}\overset{\mathrm{def}}{=}\; & \frac{\partial^2
(s-c)}{\partial y_0^2}=\frac{\partial^2 s}{\partial
y_0^2}-\frac{\partial^2 c}{\partial y_0^2}=D_{3y}-D_{2y}\\
E_1\overset{\mathrm{def}}{=}\; & \frac{\partial^2 a}{\partial
x_0\partial y_0}=-\frac{\partial}{\partial
y_0}\left(\frac{x_1-x_0}{a}\right)=-\frac{(x_1-x_0)(y_1-y_0)}{a^3}\\
E_2\overset{\mathrm{def}}{=}\; & \frac{\partial^2 c}{\partial
x_0\partial y_0}=\frac{\partial}{\partial
y_0}\left(\frac{x_0-x_2}{c}\right)=-\frac{(x_0-x_2)(y_0-y_2)}{c^3}\\
E_3\overset{\mathrm{def}}{=}\; & \frac{\partial^2 s}{\partial
x_0\partial y_0}=\frac{1}{2}\left(\frac{\partial^2 a}{\partial
x_0\partial y_0}+\frac{\partial^2 b}{\partial x_0\partial
y_0}+\frac{\partial^2 c}{\partial x_0\partial
y_0}\right)=\frac{1}{2}(E_1+E_2)\\
E_4\overset{\mathrm{def}}{=}\; & \frac{\partial^2 (s-a)}{\partial
x_0\partial y_0}=\frac{\partial^2 s}{\partial x_0\partial
y_0}-\frac{\partial^2 a}{\partial x_0\partial y_0}=E_3-E_1\\
E_5\overset{\mathrm{def}}{=}\; & \frac{\partial^2 (s-c)}{\partial
x_0\partial y_0}=\frac{\partial^2 s}{\partial x_0\partial
y_0}-\frac{\partial^2 c}{\partial x_0\partial y_0}=E_3-E_2
\end{align}
}

%********* LITERATURVERZEICHNIS ***********

%********* LITERATURVERZEICHNIS ***********

\end{document}